\documentclass{article}
\usepackage{amsmath,amssymb,amsfonts}

\newtheorem{theorem}{Theorem}

\newtheorem{corollary}{Corollary}
\newtheorem{lemma}{Lemma}

\def\gf{{\rm GF}}
\def\height{{\rm ht}}

\begin{document}

\title{Cubic Laurent Series in Characteristic 2 with Bounded Partial Quotients}
\author{David P. Robbins}
\maketitle

\begin{abstract}
There is a theory of continued fractions for Laurent series in
$x^{-1}$ with coefficients in a field $F$.  This theory bears a
close analogy with classical continued fractions for real numbers with
Laurent
series playing the role of real numbers and the sum of the terms of
non-negative degree in $x$ playing the role of the integral part.

In this paper we survey the Laurent series $u$, with coefficients in a
finite extension of $\gf(2)$, that satisfy an irreducible equation of
the form
$$  a_0(x)+ a_1(x)u+a_2(x)u^2 + a_3(x)u^3=0 $$
with $a_3 \ne 0$ and 
where the $a_i$ are polynomials of low degree in $x$ with coefficients in
$\gf(2)$.  We are particularly interested in the cases in which the
sequence of partial quotients is bounded (only finitely many distinct
partial quotients occur).  We find that there are three essentially
different cases when the $a_i(x)$ have degree $\le 1$.  We also make some
empirical observations concerning relations between different Laurent series roots
of the same cubic.

\end{abstract}

\section{Introduction}\label{sec:intro}			

Let $F$ be a field and $F(x)$ be the field of rational functions over
$F$ in an indeterminate $x$.  Let $E=F((x^{-1}))$ be the field of
formal Laurent series
\[
u=a_mx^m+a_{m-1}x^{m-1} + a_{m-2}x^{m-2}+\cdots
\]
in $x^{-1}$ with coefficients in $F$.  If $a_m \ne 0$ we say that the
degree of $u$ is $m$.  $E$ is a topological field.  Its topology is
characterized by the property that a sequence $u_n$ of formal Laurent
series converges to zero when their degrees converge to~$-\infty$.

The relationship between $F(x)$ and the extension field $E$ bears a
close analogy to that between the rational numbers and the real
numbers with polynomials in $x$ playing the role of integers.

In particular most of the basic facts about continued fractions for
real numbers have analogies for $E$.  For a Laurent series $u$ we
define its integral part to be the sum of the terms of nonnegative
degree in~$x$.

Then we define the continued fraction expansion of $u$ with an
inductive calculation as follows.  We set $u_0=u$.  Given $u_i$ we
define $p_i$ to be the integral part of $u_i$.  If $p_i \ne u_i$, we
define find $u_{i+1}$ by $1/(u_i-p_i)$ so that $u_i$ satisfies
$$  u_i=p_i+\frac{1}{u_{i+1}}.  $$
If $p_i=u_i$ we terminate the procedure.  The $p_i$'s are called the
partial quotients the $u_i$'s the complete quotients.

The sequence of $p$'s terminates exactly when $u$ is rational.  If it
terminates with $p_n$, then we have
\[u=
p_0 +
           \cfrac{1}{p_1+
            \cfrac{1}{p_2+
             \cfrac{1}{\ddots\raisebox{-.1in}{$\quad\cfrac{1}{p_n}\quad.$}
              }}}
\]

Otherwise it makes sense to write
\[
u=p_0+
\cfrac{1}{p_1+
\cfrac{1}{p_2+
\cfrac{\ \ }{\ddots
}}}
\]
Indeed if we define $c_n$ by
\[c_n=
p_0 +
           \cfrac{1}{p_1+
            \cfrac{1}{p_2+
             \cfrac{1}{\ddots\raisebox{-.1in}{$\quad\cfrac{1}{p_n},$}
              }}}
\]
$c_n$ is called the $n$th {\it convergent} to $u$ and the sequence
$c_0,c_1,\ldots$ of convergents converges in $E$ to~$u$.

More generally in the irrational case each complete quotient has the
continued fraction expansion
\[
u_n = p_n +
          \cfrac{1}{p_{n+1}+
          \cfrac{1}{p_{n+2}+
          \cfrac{\ \ }{\ddots
}}}
\]

We say that an irrational Laurent series $u$ has bounded partial
quotients if the polynomials $p_k$ are bounded in degree and we say
that a Laurent series is algebraic if it is algebraic over $F(x)$.  It
can be proved that algebraic Laurent series whose minimum polynomials
have degree 2 always have bounded partial quotients.  In fact the
sequence of partial quotients is eventually periodic (by analogy with
the theory of continued fractions for quadratic algebraic numbers).

Baum and Sweet showed in \cite{baum-sweet} that, when $F=\gf(2)$, the
cubic equation (in $y$ with coefficients in $F(x)$)
$$  x + y + xy^3= 0 $$
has a unique Laurent series solution with coefficients in $\gf(2)$ and
that this solution has bounded partial quotients.  Their proof does
not yield a descripton of what the sequence of partial quotients is.
Later Mills and Robbins succeeded in giving a complete description of
that sequence of partial quotients in \cite{mills-robbins}.  They also
provided some examples in higher characteristic.  Nevertheless it
appears that there is still very little known about the nature of
continued fractions of algebraic Laurent series.  In particular, even
though there seem to be many examples with bounded partial quotients,
for any particular example, it may be difficult or impossible to
provide a proof.  Baum and Sweet also gave some simple examples with
unbounded partial quotients.

Algebraic Laurent series with unbounded partial quotients can also be
quite complicated even when the partial quotient sequence is recognizable.  
Such Laurent series were studied by Mills and Robbins in \cite{mills-robbins} and
by Buck and Robbins in \cite{buck-robbins} and Lasjaunias in \cite{alain}.

In this paper we survey cubic Laurent series in characteristic 2.
More precisely we report on algebraic Laurent series with coefficients
in a finite field of characteristic 2 that are solutions of an
irreducible equation of the form
$$  a_0(x)+ a_1(x)y+a_2(x)y^2 + a_3(x)y^3=0 $$
where the polynomials $a_i(x)$ have coefficients in $\gf(2)$.  We
concentrate primarily on the cases where the coefficients $a_i(x)$
have degrees~$\le 1$.  However, we also make some general observations
concerning the relationships that appear to hold between different roots of the
same cubic.

\section{Algorithms}\label{sec:algorithms}			

In this section we explain how we perform the calculations in our
survey.

Let $F$ be a finite field of characteristic 2.  In most of what
follows $F$ will be the field $\gf(2)$.  Suppose that we are given
polynomials $a_0(x)$, $a_1(x)$, $a_2(x)$, $a_3(x)$ with $a_3(x) \ne 0$
in $F[x]$.  There are at most three Laurent series $u$ in $E$, with
coefficients in an algebraic extension of $F$, that satisfy the
equation
\begin{equation} \label{eq:1}
 a_0(x)+ a_1(x)u+a_2(x)u^2 + a_3(x)u^3=0\;.
\end{equation}
Using classical Newton polygon methods, we can find the beginnings of
the Laurent series solutions for any algebraic equation.  From the
initial parts of these Laurent series we can calculate the first few partial
quotients of any solution.  When a solution has bounded partial
quotients, this method requires $O(n^2)$ field operations to find $n$
partial quotients.  However, in the case of cubic equations in
characteristic 2, the method, implicit in Mills and Robbins
\cite{mills-robbins}, allows for calculation of $n$ partial quotients
in $O(n)$ time when the degrees of the partial quotients are bounded.
We review that method here.

The key to the computation is to rewrite (\ref{eq:1}) in the form
\begin{equation} \label{eq:2}
u=\frac{Q(x) u^2 + R(x)}{S(x)u^2 + T(x)} 
\end{equation}
expressing $u$ as a fractional linear transformation of $u^2$ with
coefficients that are polynomials in~$x$.

We can assume without loss of generality that no non-constant
polynomial divides all four of $Q$, $R$, $S$ and $T$.  Let $D=QT-RS$.
If $D=0$, then $u$ is a rational function of $x$.  But we are only
interested in $u$'s for which the minimum polynomial is cubic.  So we
may assume that $D\ne 0$.  We will call the degree of $D$ the {\em
height} of the cubic Laurent series $u$ and denote this quantity
by~$\height(u)$.

Suppose that $u$ is such a Laurent series with partial quotients
$p_0,p_1,p_2,\ldots$ and complete quotients $u_0,u_1,u_2,\ldots$.
Then we have
$$ u_i=p_i+1/u_{i+1}=\frac{pu_{i+1}+1}{u_{i+1}}$$ 
for all $i\ge 0$.  This shows that $u_i$ is a fractional linear
transformation of $u_{i+1}$ where the matrix that relates them has
determinant 1.  It follows that, for any $i$ and $j$, $u_i$ and $u_j$
are related by a fractional linear transformation of determinant 1
with polynomial coefficients.

In characteristic 2, since squaring is linear, it is immediate that
$u^2$ has partial quotients $p_0^2,p_1^2,p_2^2,\ldots$ and complete
quotients $u_0^2,u_1^2,u_2^2,\ldots$.  Again, for any $i$ and $j$,
$u_i^2$ and $u_j^2$ are related by a fractional linear transformation
of determinant~1.

It follows that, for any $i$ and $j$, $u_i$ is a
fractional linear transformation of $u_j^2$ with the degree of the
determinant equal to $\height(u)$.  Suppose that, for some $i$ and $j$, we
know polynomials $Q$, $R$, $S$ and $T$ such that
\begin{equation} \label{eq:2.1}
u_i = \frac{Q u_j^2 + R}{S u_j^2 +T}  
\end{equation}
There are three useful computational principles.  First, if we know
the value of $p_j$, we can deduce that
\begin{equation}\label{eq:3}
u_i = \frac{Q (p_j^2+1/u_{j+1}^2) + R}{S (p_j^2+1/u_{j+1}^2) +T}
             = \frac{ (Qp_j^2 + R) u_{j+1}^2 + Q }{ (Sp_j^2 + T) u_{j+1}^2 + S }
\end{equation}
so that we have found the matrix relating $u_i$ and $u_{j+1}$ by
performing a suitable column operation on our matrix and exchanging
columns.

Similarly, if we know the value of $p_i$, we can deduce that
\begin{equation}\label{eq:4}
         u_{i+1}=1/(p_i+u_i)=\frac{S u_j^2 +T} { (Q+p_iS) u_j^2 + (R+p_iT) } 
\end{equation}
so that we have found the matrix relating $u_{i+1}$ and $u_j$ by
performing a suitable row operation on our matrix and exchanging rows.

Finally the main computational principle is that, if we have an
equation of the form
(\ref{eq:2.1})
with known $Q$, $R$, $S$ and $T$ and the degree of $p_j$ is also known, then
we can sometimes deduce that $p_i$ is (the integral part of) the
quotient when $Q$ is divided by $S$.  Note that from (\ref{eq:2.1}) we
always have
$$   u_i-\frac{Q}{S} = \frac{Q u_j^2 + R}{S u_j^2 +T}-\frac{Q}{S}=
\frac{-D}{S(S u_j^2 +T)}.   $$ 
Since $\deg(p_j)=\deg(u_j)$ is known, we can compute 
$\deg(Su_j^2)$.  If $\deg(Su_j^2) > \deg(T)$, then we know
$\deg(Su_j^2+T)$ and therefore $\deg(S(Su_j^2+T))$.  Finally, if this 
last degree exceeds $\height(u)$, then we can conclude that $u_i$ and $Q/S$ have
the same integral part and that therefore $p_i$ is the quotient when $Q$ is
divided by $S$.  (We remark that, with sufficiently detailed knowledge
of $u_j$, it is possible that we can deduce that
$\deg(S(Su_j^2+T))>\height(u)$ without requiring that
$\deg(Su_j^2)>\deg(T)$.  But in our computations we deduce new partial 
quotients this way only when we have the
sufficient conditions that $\deg(S)+2\deg(p_j)> \deg(T)$ and
$2\deg(S)+2\deg(T)>\height(u)$.)

Once $p_i$ is known we can perform an operation of type (\ref{eq:4})
and obtain a new relation of the form (\ref{eq:2.1}) from which we may
be able to find another partial quotient, and so forth.

When we cannot deduce the value of $p_i$ this way, we may still be able to make
progress if we know sufficiently many terms of the sequence
$p_j,p_{j+1},\ldots $.  Let us assume that $j>0$ so that
$\deg(p_j)>0$. If $\deg(S) \ge \deg(T)$, then
$\deg(S(Su_j^2+T))=2(\deg(S)+\deg(u_j))$.  Moreover subsequent
operations of type $(\ref{eq:3})$ will always yield relations of the
form (\ref{eq:2.1}) in which $\deg(S) \ge \deg(T)$, where the sequence
of $S$'s have degrees increasing by at least 2.  Thus after a few
steps of this type we will be in position to compute one or more new
partial quotients $p_i$.  If we are stuck with a case in which
$\deg(S)<\deg(T)$, then a transformation of type (\ref{eq:3}) will
yield a new relation of the form (\ref{eq:2.1}) in which $\deg(T)$ is
smaller than it was before.  But there can be only finitely many steps
of this type.  Thus eventually we will arrive at the favorable case
where $\deg(S) \ge \deg(T)$ and in a few more steps we will be able to
compute a new partial quotient.

We can now see the outline of a general computational procedure.  We
start with a relation of the form (\ref{eq:2.1}) and $i=j=0$ 
and use classical methods to compute the first few
partial quotients $p_0,p_1,\ldots$.  Then, when possible, we use the
main principle above to compute a new partial quotient $p_i$ and
adjoin it to our list of known partial quotients, if it is not already
known.  (If it is already known, we have a check on our results.)  We
then use the value of $p_i$ to obtain a new relation of the form
(\ref{eq:2.1}) where $i$ has been replaced by $i+1$.  If we cannot use
the main principle, provided $p_j$ is known, we make a transformation
of the type (\ref{eq:3}) yielding a new relation where $j$ has been
replaced by $j+1$.  Thus the first type of step produces new partial
quotients and advances $i$ while the second type of step uses old
partial quotients and advances $j$.  If the partial quotient $p_j$ is
always known when needed, we can continue indefinitely.  In particular
if $i>j$ we can continue.  We have observed empirically that on
average in this algorithm $i$ increases twice and fast as $j$ so the
production rate for partial quotients is approximately twice the
consumption rate with small local variations.  However, our initial
relation of the form (\ref{eq:2.1}) has $i=j=0$, so there can be some
difficulty getting started.  Thus we use classical methods to find a
few partial quotients for an initial supply.  After the first few
steps, $i$ seems to stay reliably ahead of $j$ so production stays
reliably ahead of consumption.

The computational procedure can also be thought of as the operation of
an automaton.  From this point of view a state of the automaton is one
of the matrices relating $u_i$ and $u_j^2$ from which no deduction of
a partial quotient is immediately possible.  We think of the use of
$p_j$ as the reading of an input by the automaton, and we regard any
partial quotients $p_i,p_{i+1},\cdots $ that can be computed as new outputs
and we view the state arrived at, after computing the new partial
quotients, as the new state.  Note that this automaton is unusual in
that its inputs come from its previous outputs.

In examples with bounded partial quotients it appears that only
finitely many states occur, so we have something like a finite state
automaton.  However, we note that this is not what is usually meant by
a finite state automaton.  The reason is that we do not see every
partial quotient being read in every state.  Instead we typically find
that, for each state, there are certain partial quotients that are
never read when we are in that state.  Moreover, it is usually the
case that if one of these were read while in that state, partial
quotients never previously seen would be output, or states never
previously seen would occur.  Indeed it seems to be the avoidance of
certain combinations of states and input polynomials that makes the
partial quotient sequence bounded.  Thus the automaton description of
the partial quotient sequence removes the algebra of the problem and
makes it combinatorial.  But it does not solve the problem since there
appears to be no simple way to prove that the unseen combinations will
never occur.

The finite automaton description does, however, lead to the
possibility of even more efficient computation of the partial quotient
sequence since we can remember every combination of input polynomial
and state that occurs and what the resulting outputs are and what the
new state is.  This way the whole process can be implemented by
look-up tables.  We have not actually used this refinement in our
computations, however.

Even without the last refinement, in a typical case with bounded
partial quotients, we can find a million or so partial quotients in
just a few seconds.

There are three additional properties of algebraic elements with
bounded partial quotients that simplify our investigation.

\begin{theorem} \label{theorem:1}
Suppose that
$$ u= \frac{Qu^2+R}{Su^2+T}  $$
and that $\deg(Su^2)>\deg(T)$ and $\deg(u)\ge 0$.  If $u$ has a
partial quotient, other than the first, with degree $> \height(u) $,
then $u$ has unbounded partial quotients.
\end{theorem}

Proof: Our argument is essentially from \cite{mills-robbins}.  We
introduce the usual non-archimedean absolute value on the field of
Laurent series in which $|x|$ is set to an arbitrary real number $>1$
and $|u|=|x|^{\deg(u)}$ for any Laurent series $u$.

If the convergent $c_n$ of $u$ is $a_n/b_n$ with $a_n$ and $b_n$ relatively prime, 
then it is known \cite{baum-sweet} that
$$  |u-a_n/b_n|= \frac{1}{ |p_{n+1}| |b_n|^2},  $$
and that, conversely, if $a$ and $b$ are relatively prime polynomials with
$$  |u-a/b|= \frac{1}{|x|^k |b|^2} $$
for some positive integer $k$, then there is a non-negative integer
$n$ with $a/b=c_n$ and $k=\deg(p_{n+1})$.   We call $k$ the accuracy of
the convergent $a/b$.

In particular since we assume that 
$u$ has degree $\ge 0$, every convergent $a/b$ of $u$ has
$|a/b|=|u|$.  

Now suppose that $c=a/b$ is a convergent of $u$ of accuracy $k>\height(u)$.
Since $|a/b|=|u|$, we have
$ |Sa^2/b^2|=|Su^2|>|T|  $
so $|Sa^2|>|Tb^2|$ and $|Sa^2+Tb^2|=|Sa^2|=|Sb^2u^2|=|Sb^2u^2+Tb^2|\ne 0.$
It follows that
\begin{align*}
    \left|u- \frac{Qa^2+Rb^2}{Sa^2+Tb^2}\right| \\
    &=\left|\frac{Qu^2+R}{Su^2+T} - \frac{Qc^2+R}{Sc^2+T}\right| \\
    &=\left|\frac{(QT-RS)(u-c)^2}{(Su^2+T)(Sc^2+T)}\right|   \\
    &=\frac{1}{|x|^{2k}}\left| \frac{(QT-RS)}{ (Sb^2u^2+Tb^2)(Sa^2+Tb^2)}\right|   \\
    &=\frac{1}{|x|^{2k}}\left| \frac{(QT-RS)}{ (Sa^2+Tb^2)^2}\right| .  \\
\end{align*}
This shows that 
$ \left(Qa^2+Rb^2\right)/\left(Sa^2+Tb^2\right) $
is a convergent of accuracy $ \ge 2k-\height(u) > k$.  It follows that there are
convergents of arbitrarily large accuracy and therefore unbounded
partial quotients.

\begin{lemma}\label{lem:1}
If in the course of our algorithm for computing the continued
fraction expansion of the cubic Laurent series $u$, we find the
partial quotient $p_i$ when $i \ge j>0$, then the complete quotient
$u_i$ satisfies an equation of the form
$$ u_i= \frac{Qu_i^2+R}{Su_i^2+T}  $$
with $\deg(Su_i^2)>\deg(T)$.
\end{lemma}

Indeed since by hypothesis we can compute $p_i$, $u_i$ is related to $u_j^2$ by a
fractional linear transformation with matrix
$$
\left[
\begin{matrix}
    Q & R\\
    S & T\\
\end{matrix}
\right]
$$
in which $\deg(Su_j^2)>\deg(T) $ and $\deg(S^2 u_j^2)>\height(u)$.  
We will only be concerned with the first condition.

If $j=i$ we are done.
However, if $j<i$, we can depart from the usual computation and perform a sequence of 
steps in which we
successivley consume the partial quotients $p_j,p_{j+1},\ldots,p_{i}$.  

In the first such step we replace $S$ with
$S'=Sp_j^2+T$, and replace $T$ with $T'=S$, and replace $u_j$ by
$u_{j+1}$.  Then we have
$$ \deg(S'u_{j+1}^2)=\deg(S)+2\deg(p_j)+2\deg(p_{j+1})>\deg(S)=\deg(T'). $$
The same argument shows that subsequent steps preserve this relationship
between $S$ and $T$.  So when $p_i$ is finally consumed, we will have $u_i$
related to $u_i^2$ by a matrix with the desired property.

\begin{corollary}\label{corr:1}
Suppose that, in the course of our calculation of the
continued fraction expansion of the cubic Laurent series $u$, we find
a partial quotient $p_i$ when $i \ge j>0$, and that we find a partial
quotient $p_k$, with $k>i$ of degree exceeding $\height(u)$.  Then $u$
has unbounded partial quotients.
\end{corollary}

In practice Corollary \ref{corr:1} can be applied quite efficiently.
For most algebraic power series, the conditions of the corollary are
met for rather small $i$, $j$ and $k$, proving that the partial
quotient sequence is unbounded.

We believe that series that cannot be ruled out this way have bounded
partial quotients although it is still difficult in any individual
case to prove that this is the case.

Here is another useful principle, for which a proof is given in the
original Baum--Sweet paper~\cite{baum-sweet}.

\begin{theorem} If $u$ and $v$ are irrational Laurent series and $u$
and $v$ are related by a fractional linear transformation (of non-zero
determinant) with coefficients that are polynomials, then $u$ has
bounded partial quotients if and only if $v$ does.
\end{theorem}

\begin{corollary}
If the Laurent series $u$ satisfies an irreducible cubic equation,
with coefficients polynomial in $x$, and $u$ has bounded partial
quotients then every irrational element of the field generated by $u$
over the field of rational functions of $x$ has bounded partial
quotients.
\end{corollary}

Indeed if $v$ is in the (cubic) field generated by $u$, then 1, $u$,
$v$ and $uv$ must satisfy a linear relation with polynomial
coefficients.  We can then solve to find $v$ as a fractional linear
transformation of~$u$.

We will apply the preceding theorem only to $1/u$ and $1+u$.

Finally we have the simple principle, observed in \cite{baum-sweet}.

\begin{theorem} If $u$ is an algebraic Laurent series with bounded
partial quotients, and if, in the continued fraction for $u$, we
replace $x$ by any polynomial $p(x)$ of positive degree in $x$, then
we obtain another algebraic continued fraction with bounded partial
quotients.
\end{theorem}

We will only be interested in substituting $x+1$ for $x$.

\section{Results}\label{sec:le1}			

Here we investigate solutions of (\ref{eq:1}) when each of the
polynomials $a_i(x)$ has coeffecients in $\gf(2)$.  We concentrate
mainly on the case that the degrees of the $a_i$'s are all $\le 1$.
(We have also used our computational methods to investigate what
happens when the $a_i$ have larger degree and make a few observations
below concerning this more general situation.)  There are 256 such
equations.  However, we are interested only in polynomials which are
irreducible over the algebraic closure of $\gf(2)$.  This leaves 96
equations.

We test each of the 96 equations for Laurent series solutions with
bounded partial quotients.  In most cases we can use Corollary
\ref{corr:1} above to eliminate the solution from contention rather
quickly.  Any series for which we find $10^6$ partial quotients
without triggering the condition of Corollary \ref{corr:1} we declare
to have ``probable bounded partial quotients''.

There are 36 polynomials in our collection that have at least one
Laurent series root with probable bounded partial quotients.  However,
from the remarks above, there is a group of twelve substitutions
generated by the substitutions
$$
\begin{array}{ccc}
  x & \rightarrow & x+1 \; ;  \\
  y & \rightarrow & y+1 \; ; \\
  y & \rightarrow & 1/y \;, \\
\end{array}
$$
which preserve degrees and the property of having bounded partial
quotients.  None of the 36 polynomials is fixed under any of these
substitutions, so there are just three orbits.  Here we give a
representative of each of the three orbits.
$$
\begin{array}{llcc}
\hbox{case A}: \quad &x + y+xy^3 &=& 0 \; ;\\
\hbox{case B}: \quad &x+xy+(1+x)y^3 &=& 0 \; ; \\
\hbox{case C}: \quad &x+(1+x)y + xy^3& =& 0 \; .\\
\end{array}
$$
We give some empirical information about the partial quotient
sequences of the solutions in each of the three cases.  Each of the
equations has three Laurent series solution.  We shall see that
in each case the three solutions have closely related continued fraction
expansions.  However, there are large qualititive differences between
the three cases.

\subsection{Case A}

Case A is the previously studied Baum--Sweet cubic.  However, previous
studies considered only the solution that has coefficients in
$\gf(2)$.  The Case A equation has two other Laurent series solutions
with coefficients in $\gf(4)$.  These do not seem to have been
studied.  These two roots are equivalent in that one is mapped to the
other by the Frobenius automorphism of $\gf(4)$.  They appear also to
have bounded partial quotients.

What follows is a description of one of the solutions in $\gf(4)$.
The reader should bear in mind that we have not proved that this
description is correct although it does seem likely that the method of
\cite{mills-robbins} could be used to construct a proof.

A reasonable measure of the complexity of such a proof is the
complexity of the automaton.  We measure this by the number of
distinct pairs that occur each consisting of an input polynomial and a
state.  For this polynomial the number seems to be 36.  By contrast
the $\gf(2)$ solution is simpler and only leads to 12 pairs.

We represent $\gf(4)$ as the extension of $\gf(2)$ generated by an
element $t$ satisfying $t^2=t+1$ and we describe the solution $u$ to
the Case A equation whose leading term is the constant $t$.  There are
nine different polynomials that occur as partial quotients.  We label
these in order of appearance with the letters $a,\ldots,i$ as follows.
$$
\begin{array}{lcl}
   a &=& t  \\
   b &=& tx  \\
   c&=&  t+x  \\ 
   d&=&  (1+t) + x^2 \\
   e&=&  x \\
   f&=&  x^2  \\
   g&=& 1+t + tx \\
   h&=& (1+t)x^2  \\
   i&=&t+(1+t)x^2 \\
\end{array}
$$
Next we define some strings of polynomials.  For each non-negative
integer $n$ we define

$x_n$ to be the list of length $(8\cdot 4^n-5)/3$ alternating $h$'s
and $b$'s of the form $hbh \ldots hbh$.

$y_n$ to be the list of length $(16 \cdot 4^n-7)/3$ of the form
$efe\ldots efe$.

$u_n$ to be the list of length $(8\cdot 4^n-5)/3$ of the form
$fef\ldots fef$.

$v_n$ to be the list of length $(16\cdot 4^n-7)/3$ of the form
$bhb\ldots bhb$.

Now here is what the partial quotient sequence looks like.  The first
two partial quotients of $u$ are $a,b$. This is followed by an
infinite sequence of finite sequences $A_0,A_1,A_2,\ldots$, with $A_i$
palindromic for all $i>0$.  The $A$'s will be defined by a somewhat
complicated recursion.  We have the initial conditions
$$ \quad A_0=cdefcb, \quad A_2=ghg,\quad A_4=gibhbig  $$
and, for $n$ odd, explicit formulas: if $n=1$ mod $4$, then 
$$ A_n = eg \; x_{(n-1)/4} \; ge $$
and if $n=3$ mod 4, then
$$ A_n = cd \; y_{(n-3)/4} \; dc. $$

Here is what happens if $n$ is even and $ \ge 6$.  If $n=0$ mod 4,
$$A_n = h_n \; gi \; v_{(n-8)/4} \; ig \; r(h_n);$$
for $n=2$ mod 4, 
$$A_n = h_n \; bc \; u_{(n-6)/4}\; cb \; r(h_n),$$
where $r$ is the operator that reverses the terms in a sequence and
$h_n$ is defined below.

For $n>0$ define the palindrome $p_n$ by $$p_n=A_0\ldots A_{2n-2}
A_{2n-1} A_{2n-2}\ldots r(A_0).$$ Also set $p_0=A_3=cdefedc$ and
$p_{-1}=cfc$.

We have
$$
\begin{array}{lcl}
 h_6 &=&gibhge  \\
\end{array}
$$
If $n\ge 8$, then,
for $n=0$ mod 4,
$$   h_{n}=h_{n-2} \; bc \;  u_{(n-8 )/4} \; cb \;  p_{(n-10)/2} $$
and, for $n=2$ mod 4,
$$  h_{n}=h_{n-2} \; gi \; v_{(n-10)/4} \; ig \; p_{(n-10)/2}.  $$

This completes the recursive description of the pattern of partial
quotients.  This pattern has been verified to continue for one million
partial quotients.

The continued fraction expansion of the solution to Case A with
coefficients in $\gf(2)$ is described in \cite{mills-robbins}.  There
it is proved that the partial quotients follow a pattern somewhat
similar to the one given here.  But the connection between the
patterns is actually much more striking.  We have observed empirically
that we obtain the partial quotient sequence for the $\gf(2)$ solution
by replacing every non-zero coefficient of every partial quotient of
the $\gf(4)$ solution with a 1.  We have not found an explanation for
this phenomenon.

This phenomenon does not appear to be restricted to the Baum--Sweet
cubic.  We have observed several other cases of equations of the form
of (\ref{eq:1}) that have two roots with bounded partial quotients in
$\gf(4)$ and a root in $\gf(2)$.  In these examples we allowed the
degrees of the $a_i(x)$ to exceed 1.  In each case the root in
$\gf(2)$ also had bounded partial quotients and was related to the
$\gf(4)$ roots, but we have not been able to give a precise description of
what that relationship is.
 
\subsection{Case B}

It appears that neither Case B nor Case C has been previously studied.
They both have three solutions with coefficients in $\gf(8)$.  In each
case all three are equivalent under the Frobenius automorphism of
$\gf(8)$.

Case B is unusual in that all its partial quotients (except the first,
which is constant) have degree 1.  It also has the unusual property
that, in the finite automaton, after the first few inputs, every input
produces precisely two outputs.  In a sense this example is much more
complicated than Case A since there are 737 distinct input-state pairs
that occur.  On the other hand inspection of the list of pairs shows
that there are a great many symmetries and much structure to the list
so the complication may not be quite so great.

We can conjecture a recursion for the sequence of partial quotients.
We represent $\gf(8)$ as the field generated over $\gf(2)$ by a solution 
$t$ to $1+t+t^3=0$.
For brevity we will identify elements of $\gf(8)$ with the integers
from 0 to 7, according to their binary expansions.  Thus we will
denote $0$, $1$, $t$, $1+t$, $t^2$, $1+t^2$, $t+t^2$, $1+t+t^2$
respectively by 0, 1, 2, 3, 4, 5, 6, 7.

Also, for brevity, we will denote a polynomial $a+bx+cx^2+\cdots$ by
the sequence of digits $abc\ldots$.  So for example 13 stands for the
polynomial~$1+(1+t)x$.

Using this notation we find that the first 4 partial quotients,
$p_0,p_1,p_2,p_3$ are
$$2, \;13,\; 13,\; 01.$$
Thereafter, if we group the remaining partial quotients in quadruples,
$$ (p_4,p_5,p_6,p_7),\; (p_8,p_9,p_{10},p_{11}), \ldots, $$
there are precisely 63 quadruples that occur.

Also the 12 partial quotients $(p_4,\ldots,p_{15})$ are
$$ 33,\; 11,\; 73,\; 04,\; 53,\; 23,\; 41,\; 07,\; 11,\; 77,\; 21,\; 05 \,. $$
Thereafter, if we group the remaining partial quotients in 16-tuples,
$$  (p_{16},p_{17},\ldots,p_{31}), (p_{32},\ldots,p_{47}), \ldots, $$
there are precisely 63 16-tuples that occur.

Finally there is a bijection between the set of quadruples and the set
of 16-tuples such that, after applying the bijection the sequence of
quadruples beginning with $(p_4,p_5,p_6,p_7)$ becomes the sequence of
16-tuples, beginning with $(p_{16},\ldots,p_{31})$.  The list of 63
quadruples, together with their bijectively associated 16-tuples, is
given in Table \ref{ta:1} below.  One can use this table, together
with the initial conditions above, to generate the sequence of partial
quotients. For example, since we have
$(p_4,p_5,p_6,p_7)=(33,11,73,04)$, we can deduce that
$(p_{16},\ldots,p_{31})$ is the associated 16-tuple
$(61,03,\ldots,54,02)$.

It is not hard to find algebraic relationships in Table \ref{ta:1},
particularly if we group the rows according to the positions of zeroes
in each row.  However these algebraic relations have not yielded
additional insight and we omit them.

\begin{table}
\tiny
\caption{Case B: quadruples and 16-tuples}\label{ta:1}
\begin{verbatim}
                                11 33 41 07   54 02 64 03 44 33 34 06 14 12 76 04 66 44 26 05 
                                11 55 61 03   57 02 47 07 77 77 17 01 27 52 64 03 44 33 34 06 
                                11 77 21 05   56 02 76 04 66 44 26 05 36 62 47 07 77 77 17 01 
                                12 16 54 02   55 11 35 06 55 66 25 05 45 71 73 04 33 44 63 03 
                                13 01 73 04   53 23 41 07 11 77 21 05 71 43 34 06 44 66 54 02 
                                14 12 76 04   55 55 15 01 55 11 35 06 75 45 61 03 11 33 41 07 
                                15 01 35 06   51 27 72 04 22 44 32 06 62 37 17 01 77 11 57 02 
                                16 14 32 06   55 66 25 05 55 55 15 01 65 36 42 07 22 77 72 04 
                                17 01 57 02   52 24 63 03 33 33 13 01 43 74 26 05 66 55 56 02 
                                21 05 61 03   36 62 47 07 77 77 17 01 27 52 64 03 44 33 34 06 
                                22 44 32 06   35 06 25 05 55 55 15 01 65 36 42 07 22 77 72 04 
                                22 66 42 07   32 06 42 07 22 77 72 04 12 16 54 02 44 22 64 03 
                                22 77 72 04   34 06 54 02 44 22 64 03 74 46 25 05 55 55 15 01 
                                23 51 57 02   33 44 63 03 33 33 13 01 43 74 26 05 66 55 56 02 
                                24 53 13 01   33 33 13 01 33 11 73 04 53 23 41 07 11 77 21 05 
                                25 05 15 01   37 67 21 05 11 55 61 03 51 27 72 04 22 44 32 06 
                                26 05 56 02   31 65 56 02 66 22 76 04 46 75 15 01 55 11 35 06 
                                27 52 64 03   33 11 73 04 33 44 63 03 23 51 57 02 77 22 47 07 
                                31 65 56 02   66 55 56 02 66 22 76 04 46 75 15 01 55 11 35 06 
                                32 06 42 07   65 36 42 07 22 77 72 04 12 16 54 02 44 22 64 03 
                                33 11 73 04   61 03 41 07 11 77 21 05 71 43 34 06 44 66 54 02 
                                33 33 13 01   63 03 13 01 33 11 73 04 53 23 41 07 11 77 21 05 
                                33 44 63 03   64 03 34 06 44 66 54 02 24 53 13 01 33 11 73 04 
                                34 06 54 02   67 31 35 06 55 66 25 05 45 71 73 04 33 44 63 03 
                                35 06 25 05   62 37 17 01 77 11 57 02 37 67 21 05 11 55 61 03 
                                36 62 47 07   66 44 26 05 66 55 56 02 16 14 32 06 22 66 42 07 
                                37 67 21 05   66 22 76 04 66 44 26 05 36 62 47 07 77 77 17 01 
                                41 07 21 05   14 12 76 04 66 44 26 05 36 62 47 07 77 77 17 01 
                                42 07 72 04   12 16 54 02 44 22 64 03 74 46 25 05 55 55 15 01 
                                43 74 26 05   11 55 61 03 11 33 41 07 31 65 56 02 66 22 76 04 
                                44 22 64 03   13 01 73 04 33 44 63 03 23 51 57 02 77 22 47 07 
                                44 33 34 06   17 01 57 02 77 22 47 07 67 31 35 06 55 66 25 05 
                                44 66 54 02   15 01 35 06 55 66 25 05 45 71 73 04 33 44 63 03 
                                45 71 73 04   11 33 41 07 11 77 21 05 71 43 34 06 44 66 54 02 
                                46 75 15 01   11 77 21 05 11 55 61 03 51 27 72 04 22 44 32 06 
                                47 07 17 01   16 14 32 06 22 66 42 07 52 24 63 03 33 33 13 01 
                                51 27 72 04   44 66 54 02 44 22 64 03 74 46 25 05 55 55 15 01 
                                52 24 63 03   44 33 34 06 44 66 54 02 24 53 13 01 33 11 73 04 
                                53 23 41 07   44 22 64 03 44 33 34 06 14 12 76 04 66 44 26 05 
                                54 02 64 03   45 71 73 04 33 44 63 03 23 51 57 02 77 22 47 07 
                                55 11 35 06   42 07 72 04 22 44 32 06 62 37 17 01 77 11 57 02 
                                55 55 15 01   41 07 21 05 11 55 61 03 51 27 72 04 22 44 32 06 
                                55 66 25 05   47 07 17 01 77 11 57 02 37 67 21 05 11 55 61 03 
                                56 02 76 04   46 75 15 01 55 11 35 06 75 45 61 03 11 33 41 07 
                                57 02 47 07   43 74 26 05 66 55 56 02 16 14 32 06 22 66 42 07 
                                61 03 41 07   27 52 64 03 44 33 34 06 14 12 76 04 66 44 26 05 
                                62 37 17 01   22 44 32 06 22 66 42 07 52 24 63 03 33 33 13 01 
                                63 03 13 01   24 53 13 01 33 11 73 04 53 23 41 07 11 77 21 05 
                                64 03 34 06   23 51 57 02 77 22 47 07 67 31 35 06 55 66 25 05 
                                65 36 42 07   22 66 42 07 22 77 72 04 12 16 54 02 44 22 64 03 
                                66 22 76 04   25 05 15 01 55 11 35 06 75 45 61 03 11 33 41 07 
                                66 44 26 05   21 05 61 03 11 33 41 07 31 65 56 02 66 22 76 04 
                                66 55 56 02   26 05 56 02 66 22 76 04 46 75 15 01 55 11 35 06 
                                67 31 35 06   22 77 72 04 22 44 32 06 62 37 17 01 77 11 57 02 
                                71 43 34 06   77 11 57 02 77 22 47 07 67 31 35 06 55 66 25 05 
                                72 04 32 06   74 46 25 05 55 55 15 01 65 36 42 07 22 77 72 04 
                                73 04 63 03   71 43 34 06 44 66 54 02 24 53 13 01 33 11 73 04 
                                74 46 25 05   77 77 17 01 77 11 57 02 37 67 21 05 11 55 61 03 
                                75 45 61 03   77 22 47 07 77 77 17 01 27 52 64 03 44 33 34 06 
                                76 04 26 05   75 45 61 03 11 33 41 07 31 65 56 02 66 22 76 04 
                                77 11 57 02   73 04 63 03 33 33 13 01 43 74 26 05 66 55 56 02 
                                77 22 47 07   76 04 26 05 66 55 56 02 16 14 32 06 22 66 42 07 
                                77 77 17 01   72 04 32 06 22 66 42 07 52 24 63 03 33 33 13 01 
\end{verbatim}

\end{table}
\normalsize

\subsection{Case C}

Case C appears to have bounded partial quotients but we have not been
able to identify the pattern of partial quotients.  Other than the
first partial quotient which is a constant, the polynomials that occur
as partial quotients comprise exactly the set of all polynomials of
degree 1 together with the squares of all polynomials of degree 1.
Thus, ignoring the first partial quotient, there are 112 possible
partial quotients.

This case seems to be of much greater complexity than the others.
Over 17000 input-state pairs occur during the generation of the first
four million partial quotients and it seems as if one would have to
compute many more before all possible pairs would occur.  This casts
some doubt on the boundedness of the partial quotients.

Even though we cannot give a simple description of the partial
quotient sequence, the sequence itself is far from random looking.
For example, it contains very long subsequences which alternate
between a multiple of $x$ and a multiple of $x^2$.  These subsequences
are the centers of even larger palindromic subsequences.  The lengths
of these palindromic subsequences appear to be unbounded.

In Table \ref{ta:2}, we have listed the first 1000 partial quotients
for the solution whose constant term is $t$ in the same notation as we
used for Case B.  (The first row of the table contains the first 20
partial quotients, the second row the second twenty, etc.)

\begin{table}
\tiny
\caption{Case C: first 1000 partial quotients} \label{ta:2}
\begin{verbatim}
          2    17   52   35   77   32   36   31   73   206  13   73   74   47   05   67   21   001  01   201  
          41   03   31   66   64   43   62   74   24   53   76   46   25   002  25   46   46   25   002  05   
          702  15   15   102  75   23   04   53   76   26   57   05   47   74   73   13   306  33   36   35   
          15   702  05   002  05   702  15   35   36   13   02   43   34   007  34   43   12   76   003  06   
          703  26   44   14   007  04   307  74   15   62   16   63   77   61   23   67   62   52   36   31   
          53   12   56   44   74   02   64   507  04   007  14   44   76   55   73   36   05   16   51   57   
          004  07   004  07   004  07   004  07   004  57   51   16   05   36   73   55   76   64   53   16   
          32   47   604  07   004  07   604  47   32   16   73   05   53   13   006  13   53   45   61   001  
          01   201  41   53   12   66   16   003  16   66   22   31   001  01   501  71   66   64   53   76   
          66   72   14   75   61   201  21   47   67   34   05   74   207  54   56   62   35   07   15   13   
          57   74   03   14   51   501  21   57   504  07   004  17   37   23   41   001  41   23   47   76   
          61   71   05   51   62   305  22   35   07   05   47   74   73   13   506  03   006  03   006  03   
          006  03   006  03   006  03   006  03   006  03   006  03   006  03   506  13   73   74   47   05   
          07   35   22   305  62   51   05   71   61   16   05   26   22   21   31   07   71   74   26   003  
          06   003  06   003  06   003  06   003  26   74   71   07   31   21   22   26   15   71   71   15   
          56   603  06   003  06   603  56   15   71   21   02   71   73   006  73   71   52   64   007  04   
          407  34   43   62   04   62   23   306  03   006  03   306  23   62   04   52   51   55   37   004  
          37   55   11   46   003  06   503  36   37   01   27   34   55   35   46   36   54   507  34   13   
          51   21   02   71   73   006  03   506  13   03   43   77   61   23   67   12   71   33   75   002  
          05   402  25   06   45   77   32   36   31   53   12   46   27   12   71   13   44   15   62   16   
          63   27   104  07   004  37   55   61   401  01   001  21   67   75   15   002  15   75   07   35   
          12   54   63   206  03   006  03   206  63   54   12   65   21   53   12   56   44   54   37   52   
          35   37   204  07   004  57   51   66   64   43   62   74   24   43   62   54   33   57   74   03   
          14   31   35   002  05   002  05   002  05   002  05   002  05   002  05   002  05   002  05   002  
          05   002  05   002  05   002  05   002  05   002  05   002  05   002  05   002  05   002  05   002  
          05   002  05   002  35   31   14   03   74   57   33   54   62   43   24   74   62   43   64   66   
          51   57   004  07   204  37   35   52   37   54   44   56   12   53   51   41   64   67   07   27   
          65   302  55   32   02   02   32   65   21   73   706  03   006  03   006  03   006  03   006  03   
          006  03   006  03   006  03   006  03   006  03   706  73   21   65   32   02   02   32   55   302  
          65   27   07   67   64   71   07   41   55   47   14   06   24   65   56   003  06   003  06   003  
          06   003  06   003  56   65   24   06   14   47   55   41   47   67   24   72   76   103  06   003  
          06   103  76   72   24   47   05   67   21   001  21   67   75   15   002  05   302  35   41   27   
          54   04   04   54   27   41   55   37   004  07   004  07   004  07   004  07   004  37   55   41   
          27   54   04   04   54   67   47   01   37   36   503  06   003  06   503  36   37   01   67   45   
          66   67   004  67   66   55   43   006  03   406  53   57   74   53   76   66   72   14   75   61   
          201  21   47   67   34   05   14   32   005  02   205  72   07   12   11   47   37   43   43   37   
          67   64   31   35   002  35   31   74   26   003  06   603  56   25   502  65   27   07   67   64   
          21   56   06   76   43   61   57   404  27   23   41   001  41   23   37   17   004  07   504  57   
          21   501  51   14   03   74   57   13   15   07   35   62   56   54   207  74   05   34   67   47   
          21   201  61   75   14   72   66   76   53   64   66   71   501  01   001  31   22   66   16   003  
          16   66   12   53   41   201  01   001  61   45   53   13   006  13   53   05   73   16   32   47   
          604  07   004  07   604  47   32   16   53   64   76   55   73   36   05   16   51   57   004  07   
          004  07   004  07   004  07   004  57   51   16   05   36   73   55   76   14   47   05   47   74   
          73   13   206  73   31   36   32   77   35   52   57   77   36   503  06   003  46   11   55   37   
          004  37   55   61   701  51   14   03   74   57   13   15   07   35   62   56   54   307  24   41   
          51   03   31   46   703  66   73   05   03   25   52   57   17   304  07   004  07   004  07   004  
          07   004  07   004  07   004  07   004  07   004  07   004  07   004  07   004  07   004  07   004  
          07   004  07   004  07   004  07   004  07   004  07   004  07   004  07   004  07   004  07   004  
          07   004  07   004  07   004  07   004  07   004  07   004  07   004  07   004  07   004  07   004  
          07   004  07   004  07   004  07   004  07   004  07   004  07   004  07   004  07   304  17   57   
\end{verbatim}
\end{table}
\normalsize

\eject

\subsection{Equations with Three \protect\boldmath{$\gf(2)$} Solutions
          with Bounded Partial Quotients} 

It is quite possible for an equation of the form (\ref{eq:1}) to have
three Laurent series roots with coefficients in $\gf(2)$ each with
probable bounded partial quotients.  One of the simplest examples is
the equation
$$      1+x^2u+(1+x^2)u^2+xu^3 =0.  $$

In all cases like this that we have examined, the continued fractions
for the three roots seem to be roughly related to each other.  For
example, it appears that the set of partial quotients that occur
infinitely often is the same for all three roots.

It is also possible that an equation have three Laurent series
solutions in $\gf(2)$, just one of which has probable bounded partial
quotients, although such polynomials seem to be more rare than those
with all three roots having bounded partial quotients.  We have seen
no examples with three Laurent series roots in $\gf(2)$, precisely two
of which have probable bounded partial quotients.

\end{document}